\documentclass[final]{article}

\usepackage{amsfonts}
\usepackage{amsmath}
\usepackage{algorithm}
\usepackage{algorithmic}
\usepackage[usenames,dvipsnames]{color}
\usepackage{textcomp}
\usepackage{stmaryrd}
\usepackage{graphicx}
\usepackage{nameref}
\usepackage{multirow}
\usepackage{rotating}
\usepackage{pdflscape}
\usepackage{xfrac}
\usepackage[framed,numbered,autolinebreaks]{mcode}
\usepackage{listings}
\usepackage{color}
\usepackage{textcomp}
\usepackage{enumitem}
\usepackage{type1cm}
\usepackage{sidecap}
\usepackage{array}

\allowdisplaybreaks

\definecolor{listinggray}{gray}{0.9}
\definecolor{lbcolor}{rgb}{0.9,0.9,0.9}
\newcommand{\hiddencomment}[1]{}

\let\oldv\v

\newcommand{\vt}[1]{\mathbf{#1}}        
\newcommand{\mt}[1]{\mathbf{#1}}        

\def\BigOh{\mathcal{O}}         
\def\Real{\mathbb{R}}           
\def\transpose{^\mathrm{T}}     
\def\u{\vt{u}}                  
\def\v{\vt{v}}                  
\def\f{\vt{f}}                  
\def\A{\mt{A}}                  
\def\F{\mt{F}}                  

\def\rank{\operatorname{rank}}

\def\R{\mt{R}}                  

\def\Ao{\A_0}
\def\Ro{\R_0}   
\def\Rot{\Ro\transpose}

\def\cross_subscript{{cut}}
\def\nocross_subscript{{sub}}

\def\Omegax{{\Omega_\cross_subscript}} 
\def\Omegan{{\Omega_\nocross_subscript}}
\def\ax{{a_\cross_subscript}}                
\def\an{{a_\nocross_subscript}}
\def\Ax{\mt{A_\cross_subscript}}   
\def\An{\mt{A_\nocross_subscript}}
\def\Ax{\mt{A_\cross_subscript}} 

\def\BigOh{O}
\def\BigTheta{\Theta}

\def\hvH{[\sfrac{h}{H}]}
\def\Hvh{[\sfrac{H}{h}]}

\begin{document}

\title{The Discretely-Discontinuous Galerkin Coarse Grid for Domain Decomposition}

\author{
  Essex Edwards\thanks{(essex$|$rbridson)@cs.ubc.ca, University of British Columbia} 
  \and 
  Robert Bridson\footnotemark[1]
}

\maketitle

\begin{abstract}
We present an algebraic method for constructing a highly effective coarse grid
correction to accelerate domain decomposition. 
The coarse problem is constructed from the original matrix and 
a small set of input vectors that span a low-degree polynomial space,
but no further knowledge of meshes or continuous functionals is used.
We construct a coarse basis by partitioning the problem into subdomains
and using the restriction of each input vector to each 
subdomain as its own basis function.
This basis resembles a Discontinuous Galerkin basis on subdomain-sized elements.
Constructing the coarse problem by Galerkin projection,
we prove a high-order convergent error bound for the coarse solutions.
Used in a two-level symmetric multiplicative overlapping Schwarz preconditioner, 
the resulting conjugate gradient solver shows optimal scaling.
Convergence requires a constant number of iterations, independent of fine problem size, 
on a range of scalar and vector-valued second-order and fourth-order PDEs.
\end{abstract}



\pagestyle{myheadings}
\thispagestyle{plain}
\markboth{EDWARDS AND BRIDSON}{DDG Coarse Grid Correction}

\section{Introduction}
\label{sec:Introduction}

Many discretizations of elliptic partial 
differential equations (PDEs) lead to sparse symmetric positive definite (SPD)
linear systems of the form \(\A\u=\f\).
For large 3D problems, iterative solvers such as 
preconditioned conjugate gradient (CG) are usually necessary.
The condition number often grows quickly as $h \to\nobreak 0$,
where $h$ is the mesh element size used for discretization:
the quality of the preconditioner becomes the crucial factor for efficiency and robustness.

With an optimal preconditioner, the linear system can be solved to desired
precision in a time which scales linearly with the problem size.
Two popular and related frameworks for potentially optimal preconditioning are
multigrid (MG) \cite{briggs2000_multigrid} and domain decomposition (DD) algorithms \cite{Smith1996, Toselli2004}; 
we focus on the latter in this paper.
The key component we present is a coarse discretization of the PDE,
using larger elements of size $H$,
providing the coarse grid correction to accelerate global convergence.

A critical factor in selecting a solver is the
question of how much domain knowledge the preconditioner requires. 
Geometric approaches require the practitioner to 
re-discretize the PDE at multiple scales, which for irregular domains and/or
coefficients may be challenging.
In contrast, algebraic approaches work almost entirely with the matrix $\A$.
While algebraic methods may be more difficult to develop, they can
provide benefits in both ease of use and in handling irregular problems.

The method we propose here is essentially algebraic, but uses additional
discrete information: we ask for a small set of \emph{generating} vectors
that span the space of degree $p$ polynomials.
We construct a coarse basis by algebraically partitioning the domain into subdomains
and using the restriction of each generating vector to each 
subdomain as its own basis function.
The resulting coarse space functions are piecewise-smooth, 
with jumps at subdomain boundaries.
From this basis, we construct a coarse problem by Galerkin projection.

We derive an error bound on the solutions to the coarse problem,
and show that it is a high-order accurate convergent coarse grid approximation
for a variety of PDEs and discretizations.
Convergence requires a limited coarsening factor $\Hvh$
and sufficiently large $p$.
Combined with DD in a Krylov method,
we observe the number of required iterations decreases rapidly with $p$,
and has reduced dependence on $\Hvh$, 
e.g.,\ maintaining optimal scaling in the case $h=H^2$.

For any finite resolution of the fine problem,
our coarse bases may or may not be interpreted as discontinuous.
However, in the limit as $h\to 0$ with $H$ fixed,
they are equivalent to the bases used in 
Discontinuous Galerkin (DG) methods \cite{Arnold2002}.
We call our coarse basis functions discretely-discontinuous,
giving rise to the name Discretely-Discontinuous Galerkin (DDG) for the approach.

We provide both theoretical and numerical evidence that DDG provides a
convenient tool for easily constructing highly effective coarse
grid corrections for a wide range of problems, varying over the type
of discretization (e.g.,\ classic finite elements or finite differences), 
the domain (from Cartesian grids to adaptive unstructured meshes), and
the underlying PDE (e.g.,\ vector-valued elasticity and fourth-order biharmonic
problems).


\section{Related Work}
\label{sec:RelatedWork}

Our approach is closely related to the aggregation-based algebraic methods for constructing a coarse basis.
For a more thorough review of aggregation techniques in the MG context, 
we refer the reader to review paper by St\"uben \cite{Stüben2001_review_of_AMG}.
We review the most closely related ideas and the DD setting.

The performance of non-smoothed aggregation, like ours, depends critically on $\Hvh$.
The simplest aggregation algorithm produces a piecewise constant coarse space.
If $\Hvh = \BigOh(1)$,
then this preconditioner applied to the Poisson problem has condition number bounded independent of $h$
\cite{sala2004_aggregation, sala2005_aggregation_improved}.

For elasticity and higher-order PDEs,
a piecewise constant basis is insufficient.
Better aggregation techniques have been derived by requiring additional user input:
the vectors that span the (near-)nullspace of the PDE \cite{vanvek1996_aggregate_biharmonic},
e.g.,\ the rigid modes for elasticity and the linear polynomials for biharmonic problems.

These techniques are already optimal, in the sense that 
the number of iterations is bounded independent of problem size.
Improvements to the iteration count can come in the form of constant factor reductions
and reduced dependence on $\Hvh$ (or geometric dependencies such as the PDE, domain, coefficients, etc.).
Many works present modifications to aggregation-based techniques that improve their performance in these ways.

Despite the optimal scaling of non-smoothed aggregation, 
when $\Hvh$ is large, the aggregation-based coarse solution is a poor approximation to the actual solution.
Galerkin projection finds a solution which is optimal in energy norm,
but the near-discontinuities at subdomain boundaries dominate the energy.
One way to reduce this dependence on $\Hvh$ is to keep the aggregation basis but apply a non-Galerkin projection, 
as in over-correction methods that apply a scaling to the Galerkin solution \cite{Blaheta1988_AMG_overcorrection,Braess1995_AMG_overcorrection}.
In practice, this significantly improves the results.

Alternatively, one can work to change the basis.
Sala et al.\ \cite{sala2005_aggregation_improved} show that the subdomains used 
for aggregation can be smaller than those used for the DD smoothing step.
Following this idea, the associated term in the bound on the condition number
reduces geometrically with $H$.
This requires some additional work to come up with the extra partitions,
and enlarges the size of the coarse problem.

Another alternative is smoothed aggregation,
which smooths the basis functions, thus reducing the steep jumps at subdomain boundaries.
For the Poisson problem, 
this transforms an $\sfrac{H}{h}$ term in the condition number bound into 
an $\sfrac{H}{d}$ term, where $d$ is the smoothing diameter \cite{sala2004_aggregation}.
This keeps the size of the coarse problem the same as basic aggregation,
but requires additional work to smooth the basis,
and can increase the number of nonzeros in the coarse matrix.

Our method reduces to non-smoothed aggregation if the only generating vector is the constant vector,
and our method inherits the upper bounds proven for non-smoothed aggregation.
We increase the performance beyond non-smoothed aggregation by using a higher-order basis,
which creates to a high-order accurate rediscretization of the input PDE.
By using a $p^\mathrm{th}$ order coarse basis we reduce the energy at subdomain boundaries 
from $\sfrac{H}{h}$ to $\sfrac{H^{p+1}}{h}$.
We find the added power of higher-order bases greatly reduces the required number of iterations.



Beyond non-smooth aggregation, discontinuous functions have appeared within DD algorithms before.
For example, the restricted additive Schwarz method 
generates discontinuities at subdomain boundaries
and has improved performance relative to the equivalent smooth method \cite{Cai99arestricted,Efstathiou2003}.
When using Discontinuous Galerkin (DG) discretizations, 
discontinuities are already present in the fine problem.
Previous works have developed DD solvers specifically for DG \cite{Durta2000,Dryja2007},
or agglomorated fine DG problems to construct coarse ones \cite{Bassi2012_agglmorating_DG_elements}.
While our approach uses a basis like that of DG methods,
it does not require that the fine problem be discretized with DG,
but can be interpreted as a rediscretization of the problem 
using a DG basis on elements of size H.


\section{Preliminaries}
\label{sec:Preliminary}

The DDG algorithm requires no geometric interpretation or information for implementation
(the generating vectors which typically would contain a basis for polynomials in the nodal
coordinates are treated as a black box).
However, our understanding and analysis is intimately tied to a geometric interpretation,
so we frequently refer to the geometric properties for simplicity.
For reference, table \ref{table:symbols} lists the major symbols used throughout this paper.
They are also each defined at first use.

\begin{table}[ht]
\centering
\caption{Common symbols used throughout the paper.}
\begin{tabular}{|lcl|}
\hline
input, fine grid system  &   & $\A\u=\f$  \\
input, generating basis (tall matrix) &   & $\F$  \\
restriction to coarse space   &   &  $\Ro$ \\
restriction to overlapping subdomain $i>0$ &   & $\widetilde{\R}_i$  \\
coarse grid matrix &   & $\Ao=\Ro\A\Rot$  \\
coarse grid system &   & $\Ao\u_0=\Ro\f$  \\
element diameter in fine grid & & $h$ \\
element diameter in coarse grid & & $H$ \\
coarsening factor & & $\Hvh$ \\
polynomial degree used in coarse basis &   & $p$  \\
subdomain overlap, in geometric distance &   & $\delta$  \\
subdomain overlap, in algebraic graph distance &   & $\Delta$  \\
spatial dimension &   & $d$  \\
\hline
\end{tabular}
\label{table:symbols}
\end{table}

%

\section{The Coarse Grid}
\label{sec:CoarseGrid}
\label{sec:CoarseGridConstruction}

\begin{figure}[ht]
\centering
\includegraphics[width=0.32\textwidth]{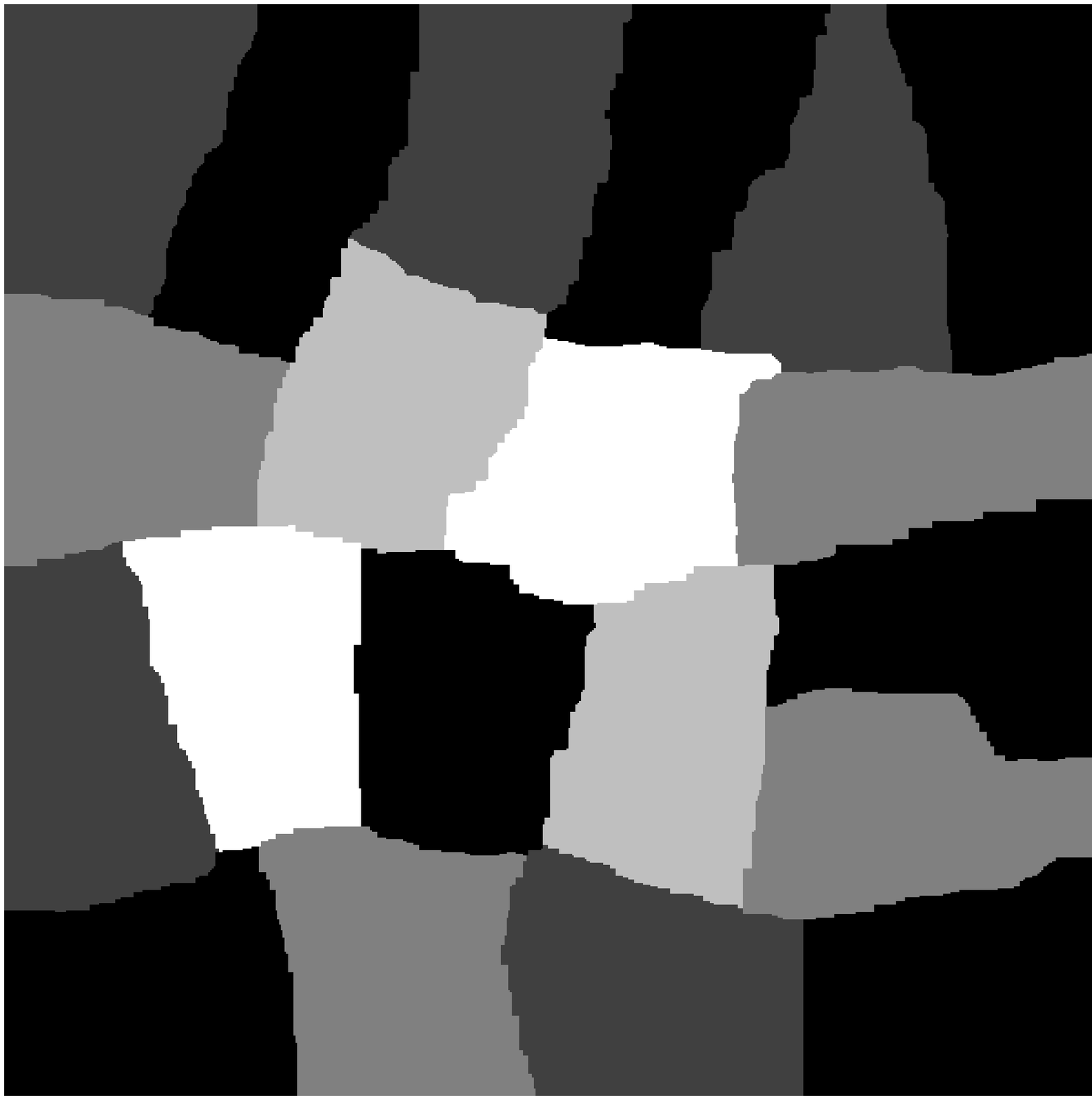}
\includegraphics[width=0.32\textwidth]{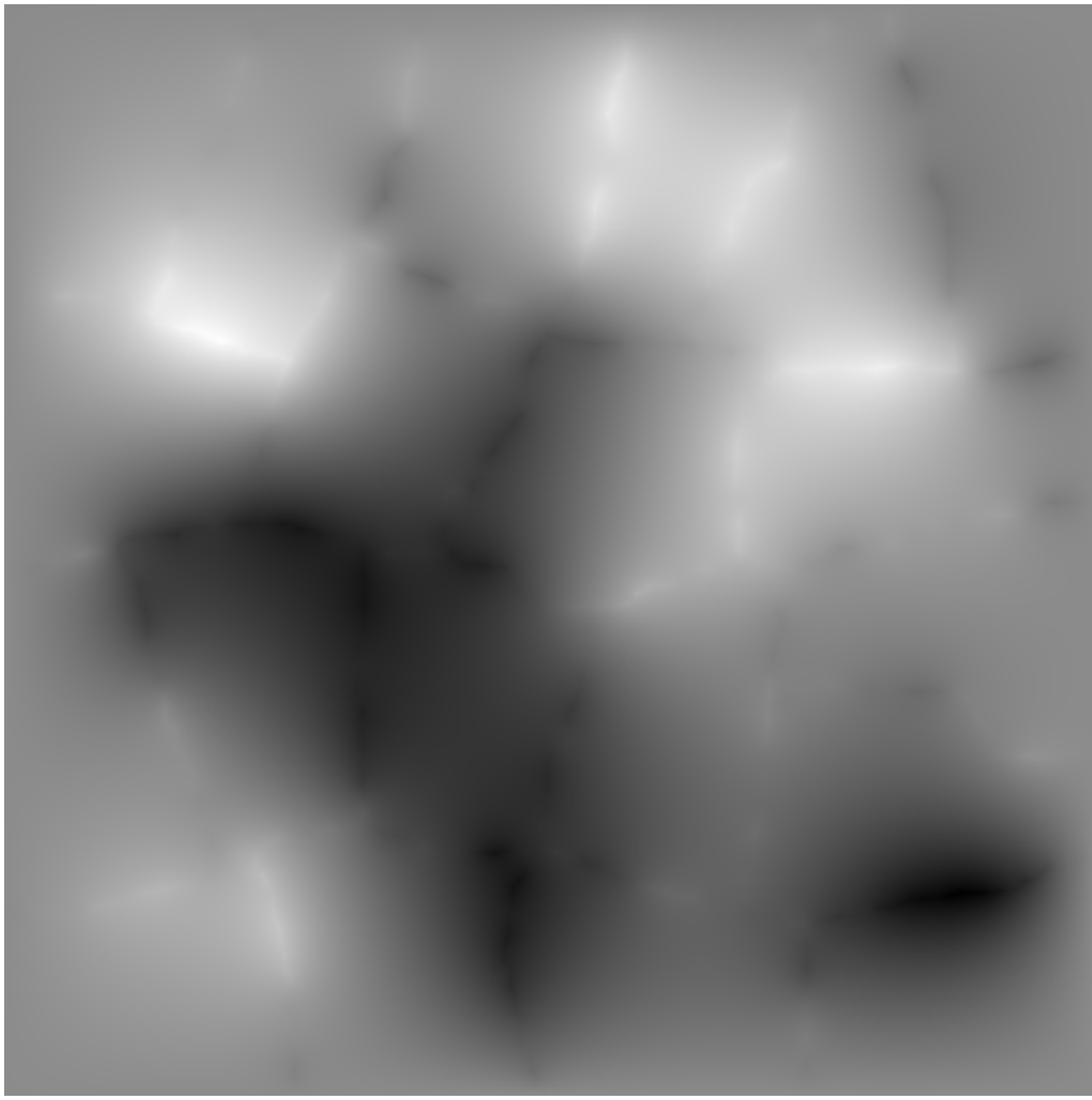}
\includegraphics[width=0.32\textwidth]{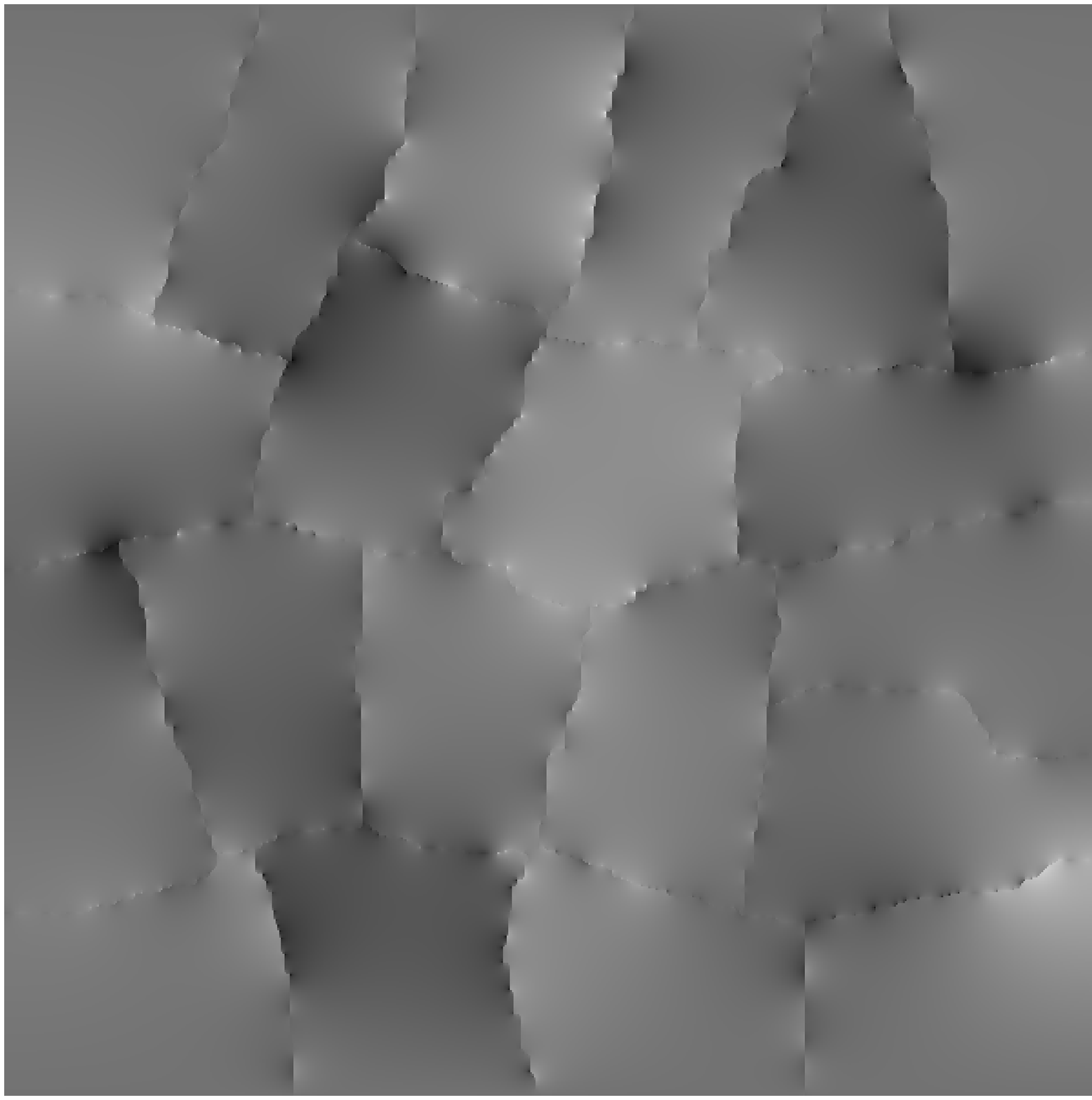}
\caption{For a simple Poisson problem with random right-hand-side, 
one-level DD produces a piecewise-smooth error after one iteration.
From left to right: partition, error after DD smoothing, $x$-derivative of error}
\label{fig:one_level_dd_error}
\end{figure}

An effective coarse grid needs to be able to approximate the error left by the one-level DD method.
After a single pass of one-level DD, the error is extremely smooth in each subdomain, 
but not across subdomain boundaries (figure \ref{fig:one_level_dd_error}).
From this structure, 
we are motivated to use piecewise higher-order polynomials which have a similar piecewise-smooth structure.

Our coarse approximation space consists of piecewise polynomials which are smooth within each subdomain,
but have arbitrary jumps at subdomain boundaries.
The construction uses several user-supplied vectors, 
arranged in the columns of $\F$, that span a degree $p$ polynomial space.
For many discretizations, $\F$ can be easily built using the nodal coordinates of the mesh and the constant vector.

To construct the coarse restriction and coarse system matrix,
we follow the same Galerkin projection as used in previous aggregation methods (e.g.\ \cite{vanvek1996_aggregate_biharmonic});
the difference is in our choice of $\F$.
The subdomains are built by partitioning the discrete domain into
non-overlapping subdomains $\Omega_i$ (i.e.\ subsets of indices) 
containing approximately $(\sfrac{H}{h})^d$ nodes for problems in dimension $d$.
Unless otherwise indicated, 
all of our examples were partitioned using 
with a graph-based discrete algorithm from the SCOTCH library \cite{Pellegrini1996_scotch}.

The coarse basis for subdomain $i$ is spanned by the columns of 
$\phi_i = \R\transpose_i \R_i\F$,
where $\R_i$ is the restriction to the 
$i^\mathrm{th}$ partition domain (with no overlap).
The final coarse restriction is made by 
concatenating these basis vectors together as rows in $\Ro$.
To help with conditioning, we orthogonalize $\Ro$,
which may be done independently for each subdomain as there is no overlap.
Orthogonalization is optional, but the remainder of the paper assumes $\Ro$ is orthogonal to simplify the analysis.
Figure \ref{fig:MATLAB_coarse_basis} shows simple (albeit inefficient) MATLAB code for this construction.
A robust implementation must also detect when $\phi_i$ is not full rank,
and discard columns as necessary.

The coarse matrix $\Ao$ is constructed by Galerkin projection, $\Ao=\Ro\A\Rot$,
and the coarse approximation to $\u$ (of $\A\u=\f$) is given by $\Rot\u_0$ where $\Ao\u_0=\Ro\f$.

\begin{figure}
\begin{lstlisting}
% F(i,j) = input basis j evaluated at node i.
% Node i is in subdomain partition(i) (zero-based index).
function R0 = basis(F,partition)
    k = size(F,2);
    [i,j,s] = find(F);
    Rt = sparse(i, partition(i)*k+j,s);
    [Rt,~]=qr(Rt,0); % optional orthogonalization.
    R0 = Rt';
\end{lstlisting}
\caption{MATLAB code to build the restriction $\Ro$ from input vectors $\F$ and a partition.}
\label{fig:MATLAB_coarse_basis}
\end{figure}

The size of $\Ao$ increases with $p$, 
both in rank $\BigTheta(p^d)$ and in the number of non-zeros $\BigTheta(p^{2d})$.
However, the block sparsity pattern of $\Ao$ is independent of $p$.
It has the same sparsity pattern as the subdomain adjacency matrix,
but with each non-zero replaced with a small dense block with size dependent on $p$.
This structure allows for efficient numerical linear algebra using dense storage and operations.
The optimal choice of $p$, in terms of total work to solve the problem,
will depend on both the problem at hand and details of the implementation, but we generally
found $p=3$, cubic polynomials, is a good default.

\section{Coarse Grid Analysis}
\label{sec:CoarseGridAnalysis}


We show that, under moderate assumptions,
the error of the coarse solution is bounded by 
\begin{equation}
\|\Rot\u_0 - \u\|_\A \leq c H^{1+p-q}(1+ \Hvh^{q-1/2} ) |u|_{\widetilde{W}_\infty^{1+p}(\Omega)}
\end{equation}
for PDEs of degree $2q$,
where $c$ is independent of $h$ and $H$,
and $u$ is the smooth interpretation of $\u$ defined in the next section.
When $\Hvh = \BigOh(1)$ and $1+p>q$, 
this error converges at high-order in $H$. A convergent coarse grid approximation
naturally allows the coarse grid correction to capture all components of the error
not handled by fine grid smoothing, leading to optimality.

We present two arguments for convergence of the coarse grid.
First, we present an argument for FEM discretizations leveraging the extensive theory surrounding FEM and Sobolev norms.
Second, we give an alternative argument that depends only on some discrete algebraic properties,
which must be shown for each particular discretization.

\subsection{Error Bound Using Geometric Properties}

Here we restrict our attention to the finite element method.

Let the domain $\Omega$ be partitioned into subdomains $\Omega_i$.
Let $u$ be in the Sobolev space $W_\infty^q(\Omega)$ 
(i.e.\ it should have bounded $q^\mathrm{th}$-order weak derivatives)
and suppose $u$ is $C^\infty$ in each subdomain.
Let $u_h \in W_\infty^q$ be a FEM interpolant of this function on some mesh.
We assume that both the subdomains and the mesh elements 
satisfy all the usual regularity assumptions for meshes
with elements of diameter $H$ and $h$ respectively.
We represent $u_h$ with a discrete vector $\u$, 
assuming a nodal basis, so $\u_i = u(x_i)$.
Furthermore, we assume that the FEM interpolant satisfies
$\|u\|_{W_2^q} \leq  c\|u_h\|_{W_2^q}$,
which is true when $\|u\|_{W_2^q} = \|u_h\|_{W_2^q} + \BigOh(h)$ and 
$h$ is smaller than some $h_0$.

Let the PDE be given as a symmetric elliptic bilinear form, $a(u,u)=\int_\Omega k[u]^2 d\Omega$,
with some linear functional $k[\cdot]$ involving up to $q^\mathrm{th}$ order derivatives.
For example, $k=\nabla$ for the Poisson problem or $k=\nabla^2$ for the biharmonic problem.
We assume continuity $a(u,u)\leq c \|u\|^2_{W_2^q}$,
where $c$ denotes an arbitrary constant independent of $h$ and $H$.
Discretized with the FEM, $\u\transpose\A\u=a(u_h,u_h)$.

\begin{figure}[ht]
\centering
\includegraphics[width=0.5\textwidth]{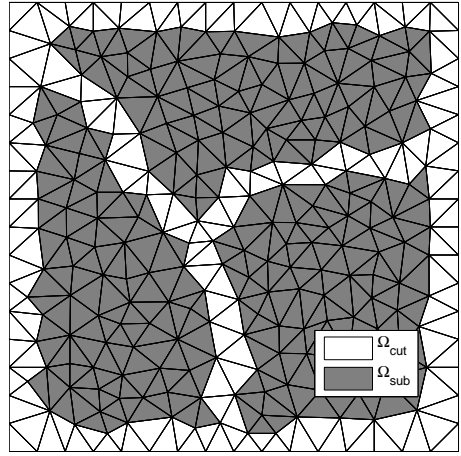}
\caption{A simple finite element mesh is divided into three subdomains.
The elements separating the subdomain interiors from each other and from the boundary
conditions are $\Omegax$.
In $\Omegan$, the error is small because the approximation is good.
The error is larger in $\Omegax$, but the total area is small enough that it does not hinder convergence.
}
\label{fig:omega_cut}
\end{figure}

Each FEM nodal point lies within exactly one subdomain, and has an associated basis function.
In some areas, basis functions from multiple subdomains overlap.
Let the union of all mesh elements containing these overlapping areas,
plus any elements touching the boundary of $\Omega$, be $\Omegax$,
and let $\Omegan = \Omega \setminus \Omegax$ (figure \ref{fig:omega_cut}).
For the purposes of the proof, 
we introduce additional bilinear forms $\ax(u,u)=\int_\Omegax k[u]^2 d\Omega$ and $\an$ defined analogously,
along with their discretizations $\Ax$ and $\An$. 
We note that $a=\ax+\an$ and $\A=\Ax+\An$.

If the mesh and subdomains are sufficiently regular
and the FEM basis functions have the usual compact support,
then the ($d$-dimensional) volume $\mu(\Omegax)\leq c\hvH$.
It is linearly dependent on $h$ because $\Omegax$ is in a band of thickness $ch$ around the subdomain boundaries,
and inversely dependent on $H$ because that is the rate at which the total subdomain surface area grows.

As in any Galerkin scheme,
the coarse solution $\Rot\u_0$ 
is the minimum error solution in the energy norm 
over the entire coarse space.
Therefore the error is bounded by that of any particular coarse vector, 
including $\v = \Rot\Ro\u$.
Using this and the splitting,
\begin{align}
\|\Rot\u_0 - \u\|_\A 
&\leq \|\v - \u\|_\A  \\
& = \left( |\v-\u|^2_\An + |\v-\u|^2_\Ax \right)^{1/2} \\
& \leq |\v-\u|_\An + |\v-\u|_\Ax
\end{align}

Before tackling either of these terms, consider what $\v=\Rot\Ro\u$ is.
Because $\Ro$ is orthogonal, $\Rot\Ro\u$ is the $l_2$ projection of $\u$ onto the coarse space.
By construction of the coarse space, $\v$ has a continuous interpretation $v$ that is 
a degree $p$ polynomial in each subdomain $\Omega_i$.
We can find $v$ directly from $u$ by a per-subdomain 
least-squares approximation of $u$ by a degree $p$ polynomial,
minimizing the sum of the squared error at each of the FEM nodal points.
Barring pathological distributions of mesh nodes, 
$v$ will be a high-order approximation to $u$,
satisfying the same error bounds commonly derived for 
FEM interpolants on a mesh with elements of size $H$.

Now, we can bound the discrete error in terms of the geometric functions.
We cite the appropriate theorems from Brenner et al.\ \cite{Brenner_2002_FEM_book}
for Sobolev and FEM-interpolant inequalities.
Looking first in $\Omegan$, 
we find that the coarse polynomials are a high-order approximation:
\begin{align}
& |\v-\u|_\An                               & \text{error interior to subdomains}   \\
& = |v_h-u_h|_{\an}                         & \text{discrete and FEM energy are equal}  \\
& \leq c \|v_h-u_h\|_{W_2^q(\Omegan)}       & \text{continuity assumption} \\
& \leq c \|v-u\|_{{W}_2^q(\Omegan)}         & \text{convergent FEM for sufficiently small $h$} \\
& = c \|v-u\|_{\widetilde{W}_2^q(\Omegan)}                   & \text{broken semi-norm defined below} \\
& \leq c \|v-u\|_{\widetilde{W}_2^q(\Omega)}                 & \text{increasing domain only increases the norm} \\
& \leq c H^{1+p-q}|u|_{\widetilde{W}_2^{1+p}(\Omega)}        & \text{Theorem 4.4.20 \cite{Brenner_2002_FEM_book}} \\
& \leq c H^{1+p-q}|u|_{\widetilde{W}_\infty^{1+p}(\Omega)}   & \text{2-norm vs. $\infty$-norm} \label{eq:bound_error_An}
\end{align}
Here the broken semi-norm $\widetilde{W}_a^b(Q)$ on domain $Q$ is defined as a sum over subdomains:
\begin{equation}
|u|_{\widetilde{W}_a^b(Q)} = \left( \sum_i |u|^a_{W_a^b(\Omega_i \cap Q)} \right)^{1/a}
\end{equation}
This is naturally extended to a maximum over subdomains for $a=\infty$.

Turning to $\Omegax$,
the coarse polynomials are not a high-order approximation in the energy norm,
because $u$ is not smooth and $v$ is not even continuous in this region.
However, $\Omegax$ is small enough that $L^\infty$ bounds are sufficient.
\begin{align}
& |\v-\u|_\Ax               & \text{error at subdomain boundaries}   \\
& = |v_h-u_h|_{\ax}         & \text{discrete and FEM energy are equal}  \\
& \leq c \|v_h-u_h\|_{W_2^q(\Omegax)}          &\text{continuity assumption} \\
& \leq c h^{-q} \|v_h-u_h\|_{L^2(\Omegax)}     &\text{Theorem 4.5.12 \cite{Brenner_2002_FEM_book}} \\
& \leq c h^{-q} \mu(\Omegax)^{1/2} \|v_h-u_h\|_{L^\infty(\Omegax)}     &\text{2-norm vs. $\infty$-norm} \\
& \leq c h^{-q} \hvH^{1/2} \|v_h-u_h\|_{L^\infty(\Omegax)}             &\text{mesh and subdomain regularity} \\
& \leq c h^{-q} \hvH^{1/2} \|v-u\|_{L^\infty(\Omegax)}                 &\text{stability of interpolation, 4.4.1 \cite{Brenner_2002_FEM_book}} \\
& \leq c h^{-q} \hvH^{1/2} H^{1+p} |u|_{\widetilde{W}_\infty^{1+p}(\Omegax)}     &\text{Theorem 4.4.20 \cite{Brenner_2002_FEM_book}} \\
& \leq c h^{-q} \hvH^{1/2} H^{1+p} |u|_{\widetilde{W}_\infty^{1+p}(\Omega)}     &\text{increasing domain} \\
& = c H^{1+p-q} \Hvh^{q-1/2} |u|_{\widetilde{W}_\infty^{1+p}(\Omega)}     &\text{factor to match eq. (\ref{eq:bound_error_An})}  \label{eq:bound_error_Ax}
\end{align}

Combining the two bounds \ref{eq:bound_error_An} and \ref{eq:bound_error_Ax} completes the error bound
\begin{align}
\|\Rot\u_0 - \u\|_\A
\leq c H^{1+p-q}(1+ \Hvh^{q-1/2} ) |u|_{\widetilde{W}_\infty^{1+p}(\Omega)}
\end{align}

\subsection{Error Bound Using Algebraic Properties}

\def\Bx{{\mt{B}_\cross_subscript}}
\def\Ix{{\mt{I}_\cross_subscript}}

We can derive a similar bound based purely on algebraic components.
As before, let $\v=\Rot\Ro\u$
and $\A=\Ax+\An$ be a splitting
into symmetric positive semi-definite components.
This splitting need not correspond to the FEM definition given earlier.
However, we require that
$\Ax = 
[\begin{smallmatrix}
  \Bx & \mt{0} \\
  \mt{0} & \mt{0}
 \end{smallmatrix}]
$
with $\Bx \in \Real^{m\times m}$
and $m\leq c \hvH h^{-d}$.
This is usually true and plays the role of $\mu(\Omegax)$ from the geometric proof.

We assume that for all $\u \in V$,
\begin{align}
& \|\Ax\|_2 \leq c h^{-2q},   \label{AssumeAxNorm} \\
& \|\u\|_\infty \leq ch^{d/2} \|\u\|_2, \\
& |\u-\v|_\An \leq cH^{1+p-q}\|\u\|_2,        \label{AssumeAn} \\
\textrm{and\;}& \|\u-\v\|_\infty \leq cH^{1+p}\|\u\|_\infty.     \label{AssumeLinf}
\end{align}
The constants $c$ include the roughness term $|u|_{\widetilde{W}_\infty^{1+p}(\Omega)}$ from the geometric proof. 
Geometrically speaking, the subspace $V$ must be restricted to functions with bounded roughness.

Showing that these assumptions are true for a particular discretization
could exploit geometric properties as in the previous section.

From these assumptions, 
the convergence argument follows the exact same structure as in the geometric case
and we do not repeat it.
The final error bound is similar to the above:
\begin{equation}
\|\Rot\u_0 - \u\|_\A \leq c H^{1+p-q}(1+ \Hvh^{q-1/2} ) \|\u\|_2.
\end{equation}

\hiddencomment{
Letting $\Ix=[\mt{I}, \mt{0}]$,
we can write $|\vt{x}|_\Ax = |\Ix\vt{x}|_\Bx$.
As before,
\begin{align}
\|\v - \u\|_A \leq |\v-\u|_\An + |\v-\u|_\Ax
\end{align}
We have simply assumed a bound on the left term by hypotheses.
To handle the right term:
\begin{align}
& |\v-\u|_\Ax & \text{} \\
& = |\Ix(\v-\u)|_{\mt{B}_\cross_subscript} & \text{by assumption} \\
& \leq \|\Bx\|_2^{1/2} \|\Ix(\v-\u)\|_2 & \text{def'n of matrix norm} \\
& \leq c h^{-q} \|\Ix(\v-\u)\|_2 & \text{assumption} \\
& \leq c h^{-q} \sqrt{m} \|\v-\u\|_\infty & \text{inf vs 2 norm} \\
& \leq c h^{-q} \hvH^{1/2} h^{-d/2} \|\v-\u\|_\infty & \text{assumption} \\
& \leq c h^{-q} \hvH^{1/2} h^{-d/2} H^{1+p} \|\u\|_\infty & \text{assumption} \\
& \leq c h^{-q} \hvH^{1/2} h^{-d/2} H^{1+p} h^{d/2} \|\u\|_2 & \text{assumption} \\
& = c h^{-q} \hvH^{1/2}  H^{1+p}\|\u\|_2 & \text{}
\end{align}
}

\subsection{Proof vs.\ Practice}

The proof and our use of the coarse grid in practice are not entirely consistent with each other.
The coarse grid is used to approximate the error after applying one-level DD.
For PDEs with smooth coefficients, as $h\to 0$ with $H$ fixed, 
the error after one-level DD is piecewise $C^\infty$ as in the proof.
For any finite $h$, it is only an approximation as accurate as the discretization.

In problems with discontinuous coefficients,
even as $h\to 0$, the error after one-level DD is not $C^\infty$ in each subdomain.
It has kinks where the coefficients have discontinuities.
We tried matching the partition boundaries to the discontinuities,
or using a piecewise generating basis $\F$ that matches the discontinuities.
We observed optimal scaling even without these strategies,
but either strategy significantly accelerated convergence with high-order polynomials.

When we use non-trivial overlap between subdomains, the boundaries of the smooth regions
do not line up with the discontinuities in the coarse space.
This is easy to resolve by adjusting the coarse subdomains,
but this introduces more variation in the coarse subdomains' size and shape.
In numerical experiments, better results were obtained by ignoring this inconsistency with the proof
and keeping the original subdomain shapes.

The roughness term in the error bound can increase with $p$,
suggesting that the error can actually increase with $p$.
However, because increasing $p$ always grows the coarse vector space,
and the Galerkin solve is optimal in that space,
increasing $p$ never increases the error.

\section{Domain Decomposition}
\label{sec:DomainDecomposition}

We combine the coarse grid within a 
standard multiplicative algebraic DD framework.

The DD subdomains begin with the same partition computed
for the construction of the coarse basis.
From the partition, overlapping subdomains are algebraically constructed 
by expanding the partition to include nodes within graph distance $\Delta$ in the graph defined by $\A$.
The expanded subdomains $\widetilde{\Omega}_i$ overlap in geometric bands of size $\delta$.
For simple meshes and discretizations, $\delta=(1+2\Delta)h$.
Unless otherwise indicated, we use minimal overlap $\Delta=0$.

Let $\widetilde{\R}_i$ be the restriction matrix, 
such that $\u_i = \widetilde{\R}_i \u$ is the vector of the elements of $\u$ from subdomain $\widetilde{\Omega}_i$,
i.e.\ $\widetilde{\R}_i$ is a subset of the rows of the identity matrix.
For each subdomain, the local problem uses the matrix $\A_i = \widetilde{\R}_i \A \widetilde{\R}\transpose_i$,
and we solve these subdomain problems exactly.
Given a current approximation $\u_k$,
processing subdomain $i$ updates the approximation to
\begin{equation}
\u_{k+1} = \u_k + \widetilde{\R}\transpose_i\A_i^{-1}\widetilde{\R}_i(\f-\A\u_k).
\end{equation}
Iterating over all subdomains and updating the approximation to $\u$ after each, 
we arrive at the algorithm for one-level multiplicative overlapping Schwarz.

To build a two-level method, 
we multiplicatively combine one pass of one-level Schwarz as a pre-smoother,
the coarse problem solution, and another pass through the subdomains as a post-smoother.
The post-smoother is done in reverse order, 
making the entire operation symmetric and usable with CG.

We give some experimental results with a three-level method operating in a V-cycle.
We construct the three-level problem by taking the two-level
algorithm and applying it again to the coarse matrix $\Ao$ to make an even coarser matrix $\A_1$.
To do this, we need a coarsened version of the generating vectors,
which are simply $\F_0 = \Ro \F$.
The coarsened coarse problem is equivalent to directly coarsening the original problem with larger subdomains.
We keep the ratio between physical element sizes in adjacent levels the same 
(i.e.\ $h/H_0 = H_0/H_1$).
The algebraic overlap $\Delta$ used for smoothing is also the same at all levels.

\subsection{Condition Number}
For several special cases, our approach reduces to previously published aggregation approaches. 
For the Poisson problem using $p=0$,
Sala \cite{sala2004_aggregation} showed that the additive variant of the preconditioner has condition 
number bounded by $\BigOh(1+\Hvh)$.
For elasticity and biharmonic problems with $p=1$, our approach is essentially a non-smoothed two-level variant of the 
multigrid method described by Van\oldv{e}k et al.\ \cite{vanvek1996_aggregate_biharmonic}.
Their coarse grid uses the zero-energy modes, which are $p=1$ for biharmonic and a subset of $p=1$ for elasticity.
They later prove optimal convergence, but only for the Poisson problem \cite{vanek2001_convergence}.
We do not have a condition number bound showing the dependence on $p$ and $q$.
However, increasing $p$ beyond the low-order choices in the literature increases the dimension of the coarse space, 
which does not have a negative effect on convergence -- it can only increase the rate of convergence.

\section{Numerical Experiments}
\label{sec:NumericalExperiments}

\begin{figure}[h]
\begin{minipage}{\textwidth}
  \centering
  $\vcenter{\hbox{\includegraphics[height=2.8in]{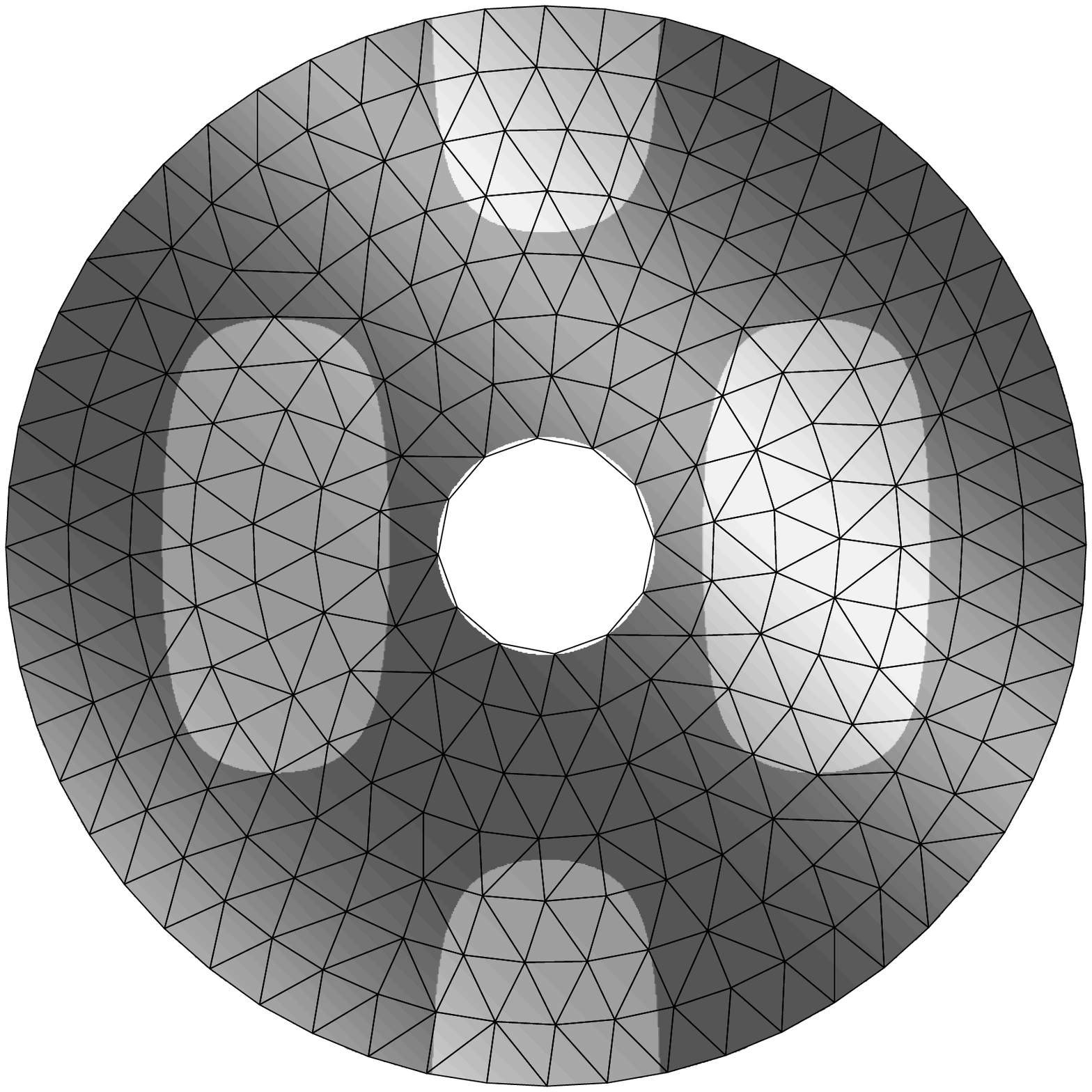}}}$
  $\vcenter{\hbox{\includegraphics[ width=2.8in, angle=-90]{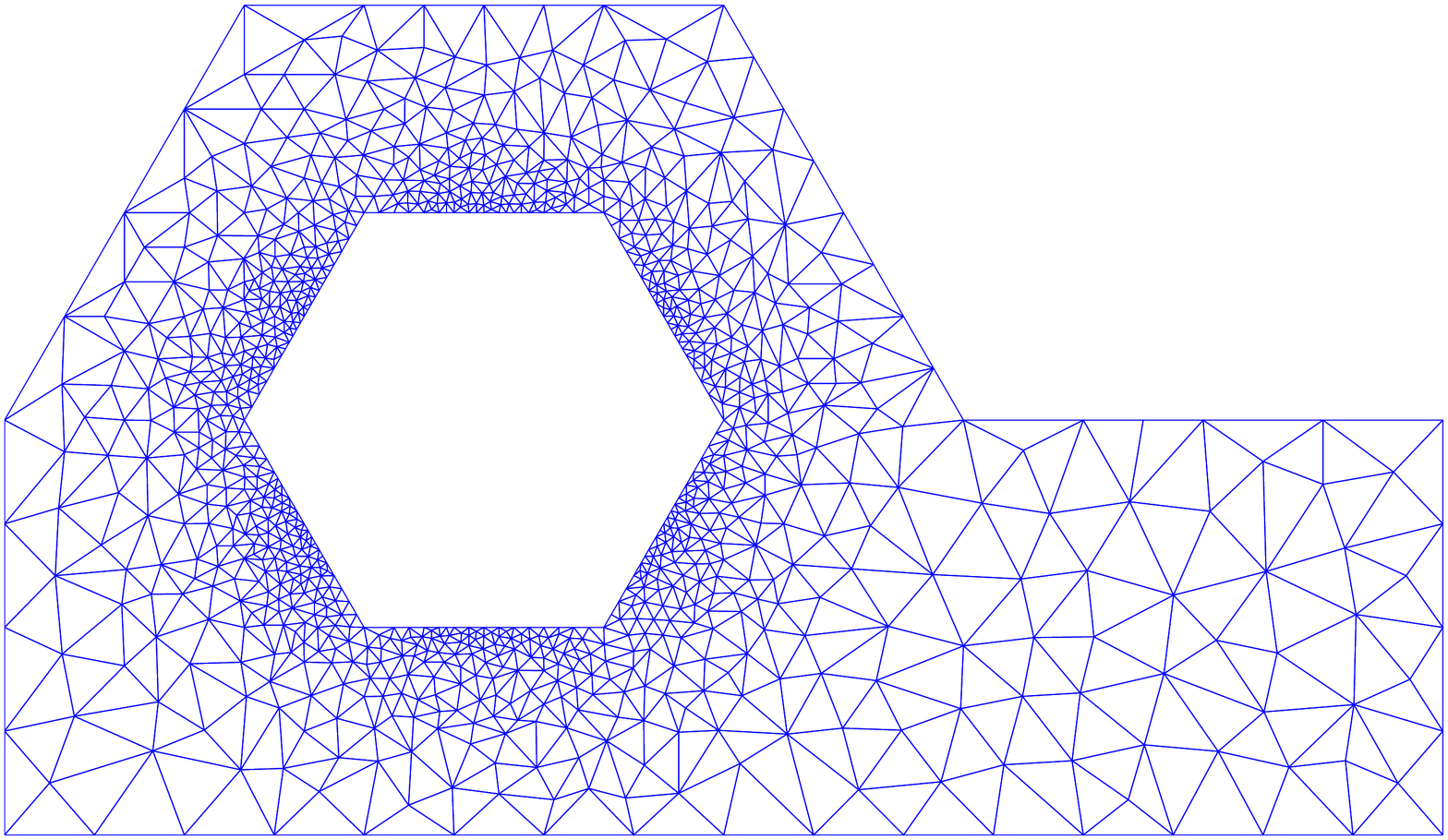}}}$
\end{minipage}
\caption{
Left: sample annulus mesh used for Poisson problems B,C,D, with shading showing $S + 10J$.
Right: sample adaptive mesh used for elasticity problem E.
}
\label{fig:FEM_meshes}
\end{figure}

We demonstrate the performance of our coarse grid and DD as a preconditioner for CG on a variety of PDEs and discretizations.
Unless otherwise noted, all problems are solved to a $10^{-9}$ reduction in 
residual after the first application of the preconditioner.
The right-hand-side vector $\f$ is a random Gaussian-distributed vector, and the initial guess for $\u$ is $\vt{0}$.
For the sake of easier reporting, 
we consider all problems with uniform meshes to be scaled such that $\operatorname{rank}(\A)=n=h^{-d}$,
where $d$ is the spatial dimension.
Consequently, $H^{-d}$ is the number of subdomains.

The graph of the logarithm of the residual vs.\ the iteration count is typically very straight.
Therefore we measure not just the integer iteration on which the residual is first smaller than the tolerance,
but also the fractional iteration count at which the linear interpolation of this graph meets the tolerance. 
We found this reveals a lot of otherwise hidden detail, and this is shown in some figures.
When CG takes more than 1000 iterations, 
we stop the solve and report iteration bounds based on condition number estimates 
derived from the Lanczos coefficients computed during CG.
For converged problems, 
these bounds agreed very well with actual iteration counts.

Scaling with $h$ and $H$ are already well explored in the existing aggregation-based literature.
Our approach does not perform significantly differently along these axes,
so we concentrate on the dependence on $p$ and the novel scaling regimes that our approach can handle.


We use the following problems:
\begin{enumerate}[label=(\Alph*)]
\item \emph{Poisson 3D.} $\nabla\cdot\nabla u = f$ discretized on an $m \times m \times m$ regular grid with a 7-point finite difference stencil.
One face of the cube has a Dirichlet boundary condition and the remainder are Neumann.
For this problem, we partition using recursive inertial partitioning \cite{Taylor94_RIP_spatial_partition} 
so that the matrix need never be explicitly constructed. 
The first partition uses a randomly-oriented plane to ensure irregularly shaped subdomains.

\item \emph{Smooth Poisson.} $\nabla\cdot S \nabla u = f$ discretized with piecewise linear finite elements on a 2D unstructured triangle mesh
of a circular annulus with outer radius 5 times the inner radius. Both inner and outer boundaries use Dirichlet conditions.
The scalar function $S=\exp( 1 + \sin(\pi(x+y)))$ is smooth. See figure \ref{fig:FEM_meshes}.

\item \emph{Non-Smooth Poisson.} $\nabla\cdot D \nabla u = f$ discretized as above, but with discontinuous 
$D=S + 100 J$ where $J=\operatorname{H}\left(0.25+\cos(\pi x)\cos(2\pi y)\right)$
with Heaviside step function $H$.
$J$ is an indicator function for two `materials' in the problem.
We used algebraic partitions that do not conform to the material boundaries,
but we use generating vectors $\F$ that are piecewise polynomial with respect to the material domains.
This doubles the number of columns in $\F$, 
but only subdomains that include the material boundary end up with additional coarse basis functions,
so $\Ao$ is only marginally larger.
In practice, we observe optimal scaling even without this extra work
and it makes no difference with the piecewise constant basis.
However, with the piecewise cubic basis, this material-aware $\F$ reduces the iteration count by nearly one half.

\item \emph{High-Order Poisson.} $\nabla\cdot\nabla$ discretized as in B and C, but with continuous piecewise cubic finite elements.

\item \emph{Elasticity.} $(\nabla\cdot\nabla + \nabla\nabla\cdot) u = f$ with Dirichlet boundary conditions, with vector $u$.
This is discretized on a spatially-adaptive unstructured 2D triangle mesh (figure \ref{fig:FEM_meshes}) with piecewise linear FEM.
For the generating vectors $\F$, we take the degree $p$ polynomials in each component of $u$.

\item \emph{Biharmonic.} $\nabla^4 u = f$ on a regular $m\times m$ 2D grid discretized with a 13 point finite difference stencil. 
All boundaries have homogenous Dirichlet and Neumann conditions.
This problem uses an algebraic overlap $\Delta=1$, since performance is quite poor with $\Delta=0$.
Also this problem is solved only to $10^{-6}$ reduction in residual,
as the fine discretizations are very poorly conditioned causing CG to break down before reaching $10^{-9}$ as used in the above.

\end{enumerate}

Table \ref{table:BigResultsTable} summarizes the results for solving these problems at different $h$, $H$, 
and polynomial coarse spaces from piecewise constant $P_0$ to cubic $P_3$.
With fixed $\Hvh$, we observe near constant iteration counts, independent of the problem size, when using the two-level algorithm.
Part of the increase in iteration count can be attributed to degrading partition quality with the algebraic partitioner.
Experiments (not shown) with more structured partitioning of the structured meshes showed less variation in iteration counts.
The three-level V-cycle does not perform nearly as well, but still appears to be sub-logarithmic in $n$.

\begin{table}[!htb]
\def\TooSmall{$\times$} 
\def\TS{\TooSmall} 
\newcolumntype{H}{>{\setbox0=\hbox\bgroup}c<{\egroup}@{}} 
\newcommand\T{\rule{0pt}{2.6ex}}       
\newcommand\B{\rule[-1.2ex]{0pt}{0pt}} 

\centering
\caption{
Iterations of CG to solve various problems in $n$ variables.
$P_0$, i.e.\ non-smoothed aggregation, is not expected to achieve optimal scaling for $\nabla^4$;
the $P_0$ column is included for comparison only.
Problems marked ``\TooSmall'' are too small to use the three-level solver - 
the coarsest partition would be a single subdomain.
}

\begin{tabular}{|rr|rrrr|rrrr|Hrrr|}
\hline
&  & \multicolumn{8}{c|}{Two-Level} & \multicolumn{4}{c|}{Three-Level} \\ 
&  & \multicolumn{4}{c|}{$H/h$=10} & \multicolumn{4}{c|}{$H/h$=20} & \multicolumn{4}{c|}{$H/h$=10} \\ 
& $n^{1/d}$ & $P_0$ & $P_1$ & $P_2$ &$ P_3$ & $P_0$ & $P_1$ & $P_2$ &$ P_3$ & $P_0$ & $P_1$ & $P_2$ &$ P_3$ \\ 
\hline 
\multirow{5}{*}{\begin{sideways} $\nabla\cdot\nabla$ 3D \end{sideways}}
\T &   40 & 36 & 20 &  15 &  12 & 35 & 23 &  18 &  15 &  36 &  \TooSmall &  \TooSmall &  \TooSmall \\
   &   80 & 41 & 20 &  16 &  13 & 51 & 28 &  21 &  18 &  41 &  \TooSmall &  \TooSmall &  \TooSmall \\
   &  160 & 44 & 21 &  16 &  13 & 61 & 30 &  23 &  19 &  73 &  39 &  29 &  23 \\
   &  320 & 44 & 21 &  16 &  14 & 63 & 30 &  23 &  19 &  92 &  43 &  32 &  26 \\
\B &  640 & 46 & 22 &  16 &  14 & 65 & 31 &  23 &  19 & 100 &  48 &  35 &  34 \\
\hline
\multirow{5}{*}{\begin{sideways} $\nabla\cdot S \nabla$ \end{sideways}}
\T &  200 & 41 & 19 &  14 &  11 & 51 & 25 &  20 &  15 &  46 &  26 &  23 &  19 \\ 
   &  400 & 44 & 19 &  14 &  11 & 60 & 26 &  19 &  16 &  69 &  40 &  29 &  24 \\ 
   &  800 & 48 & 20 &  15 &  12 & 64 & 27 &  20 &  17 &  95 &  43 &  31 &  23 \\ 
   & 1600 & 52 & 20 &  15 &  12 & 68 & 28 &  22 &  17 & 115 &  45 &  32 &  25 \\ 
\B & 3200 & 52 & 21 &  16 &  13 & 71 & 29 &  22 &  17 & 133 &  45 &  32 &  25 \\ 
\hline
\multirow{5}{*}{\begin{sideways} $\nabla\cdot D \nabla$ \end{sideways}} 
\T &   200&    43&    19&    15&    12    &    56&    27&    22&    16    &    65&    32&    27&    25 \\
   &   400&    46&    20&    15&    13    &    64&    28&    22&    19    &    71&    42&    28&    24 \\
   &   800&    48&    21&    16&    13    &    67&    31&    23&    18    &   101&    43&    31&    25 \\
   &  1600&    52&    21&    17&    14    &    69&    33&    23&    20    &   118&    47&    34&    26 \\
\B &  3200&    58&    24&    17&    13    &    69&    31&    24&    19    &     -&    49&    35&    27 \\
\hline
\multirow{5}{*}{\begin{sideways} $\nabla^2$ $P_3$-FEM \end{sideways}}
\T &  200 & 47 & 22 &  17 &  14 & 58 & 29 &  21 &  17 &  54 &  33 &  24 &  21 \\ 
   &  400 & 53 & 22 &  17 &  14 & 68 & 29 &  22 &  19 &  77 &  44 &  34 &  25 \\ 
   &  800 & 53 & 24 &  17 &  14 & 73 & 31 &  23 &  19 & 108 &  47 &  33 &  26 \\ 
   & 1600 & 56 & 23 &  19 &  15 & 75 & 32 &  25 &  19 & 133 &  49 &  35 &  28 \\ 
\B & 3200 & 67 & 24 &  18 &  15 & 83 & 32 &  24 &  20 & 147 &  50 &  36 &  28 \\ 
\hline
\multirow{5}{*}{\begin{sideways} Elasticity \end{sideways}}
\T &  200 &    72 &    29 &    21 &    19 &    83 &    40 &    32 &    25 &    66 &   \TS &   \TS &   \TS \\
   &  400 &    81 &    29 &    24 &    18 &    96 &    41 &    30 &    25 &   111 &    47 &    37 &    30 \\
   &  800 &    82 &    33 &    23 &    20 &   104 &    42 &    31 &    26 &   154 &    67 &    48 &    39 \\
   & 1600 &    76 &    32 &    24 &    20 &   107 &    44 &    33 &    27 &   195 &    69 &    51 &    36 \\
\B & 3200 &    76 &    33 &    27 &    23 &   108 &    47 &    36 &    29 &   206 &    71 &    50 &    40 \\
\hline
\multirow{4}{*}{\begin{sideways} $\nabla^4$ \end{sideways}}
\T &  200 &  698 & 62 &  20 &  12 &   600 & 154 &  44 &  24 &   722  &  119 &  40 &  22 \\ 
   &  400 & 2900 & 68 &  21 &  12 &  3500 & 188 &  53 &  27 &  4300  &  376 & 112 &  37 \\ 
   &  800 & 5700 & 77 &  25 &  15 &  6900 & 184 &  55 &  32 &  7900  &  961 & 156 &  51 \\ 
\B & 1600 & 9500 & 99 &  31 &  15 & 11000 & 246 &  70 &  30 & 14000  & 2100 & 169 &  58 \\ 
\hline
\end{tabular}
\label{table:BigResultsTable}
\end{table}

\begin{table}[!htb]
\def\NoConverge{$\times$} 
\def\ZERO{\hphantom{0}} 
\newcolumntype{H}{>{\setbox0=\hbox\bgroup}c<{\egroup}@{}} 
\newcommand\T{\rule{0pt}{2.6ex}}       
\newcommand\B{\rule[-1.2ex]{0pt}{0pt}} 

\centering
\caption{
Execution times for the middle column of table \ref{table:BigResultsTable} (two-level scheme with $H/h$=20).
Each entry shows the wall-clock time in seconds and the percent of total runtime spent on coarse grid setup and solution.
Problems marked ``\NoConverge'' did not converge within 1000 iterations.
}

\begin{tabular}{|rr|rrrr|rrrr|Hrrr|}
\hline
& $n^{1/d}$ & $P_0$ & $P_1$ & $P_2$ & $P_3$ \\ 
\hline 
\multirow{5}{*}{\begin{sideways} $\nabla\cdot\nabla$ 3D \end{sideways}}
\T &    40 &      79 ( 0\%) &      55 (\ZERO 0\%) &      42 (\ZERO 0\%) &      32 (\ZERO 1\%) \\
   &    80 &     152 ( 1\%) &      87 (\ZERO 1\%) &      64 (\ZERO 2\%) &      53 (\ZERO 3\%) \\
   &   160 &     481 ( 2\%) &     238 (\ZERO 4\%) &     186 (\ZERO 6\%) &     160 (11\%) \\
   &   320 &    1207 ( 3\%) &     614 (\ZERO 8\%) &     525 (15\%) &     497 (27\%) \\
\B &   640 &    7053 ( 5\%) &    3785 (16\%) &    3140 (25\%) &    3563 (43\%) \\
\hline
\multirow{5}{*}{\begin{sideways} $\nabla\cdot S \nabla$ \end{sideways}} 
\T &  200 &       12 ( 2\%)&        5 ( 2\%)&        4 ( 3\%)&        3 ( 7\%)\\
   &  400 &       48 ( 1\%)&       22 ( 2\%)&       17 ( 4\%)&       15 ( 8\%)\\
   &  800 &      246 ( 1\%)&      108 ( 2\%)&       87 ( 4\%)&       72 ( 9\%)\\
   & 1600 &     1283 ( 1\%)&      572 ( 2\%)&      423 ( 4\%)&      358 ( 9\%)\\
\B & 3200 &     6491 ( 0\%)&     2767 ( 2\%)&     2107 ( 4\%)&     1711 ( 9\%)\\
\hline
\multirow{5}{*}{\begin{sideways} $\nabla\cdot D \nabla$ \end{sideways}} 
\T &   200 &       11 ( 1\%)&        5 ( 2\%)&        4 ( 5\%)&        4 (11\%)\\
   &   400 &       51 ( 1\%)&       22 ( 2\%)&       18 ( 5\%)&       16 (11\%)\\
   &   800 &      227 ( 1\%)&      100 ( 2\%)&       84 ( 6\%)&       67 (13\%)\\
   &  1600 &     1051 ( 1\%)&      496 ( 2\%)&      362 ( 6\%)&      320 (12\%)\\
\B &  3200 &     7170 ( 0\%)&     3006 ( 2\%)&     2254 ( 4\%)&     1896 (10\%)\\
\hline
\multirow{5}{*}{\begin{sideways} $\nabla^2$ $P_3$-FEM \end{sideways}}
\T &   200 &       18 ( 1\%)&        9 ( 1\%)&        6 ( 4\%)&        6 ( 5\%)\\
   &   400 &       92 ( 0\%)&       43 ( 1\%)&       32 ( 3\%)&       27 ( 6\%)\\
   &   800 &      440 ( 0\%)&      190 ( 1\%)&      149 ( 3\%)&      121 ( 7\%)\\
   &  1600 &     1902 ( 0\%)&      853 ( 1\%)&      635 ( 4\%)&      498 ( 8\%)\\
\B &  3200 &    10223 ( 0\%)&     4123 ( 1\%)&     3183 ( 3\%)&     2706 ( 7\%)\\
\hline
\multirow{5}{*}{\begin{sideways} Elasticity \end{sideways}}
\T &   200 &       30 ( 0\%)&       15 ( 1\%)&       11 ( 3\%)&        9 ( 6\%)\\
   &   400 &      147 ( 0\%)&       60 ( 1\%)&       46 ( 3\%)&       42 ( 7\%)\\
   &   800 &      644 ( 0\%)&      263 ( 1\%)&      207 ( 3\%)&      172 ( 7\%)\\
   &  1600 &     2778 ( 0\%)&     1160 ( 1\%)&      883 ( 4\%)&      727 ( 8\%)\\
\B &  3200 &    14958 ( 0\%)&     6202 ( 1\%)&     5017 ( 3\%)&     4625 ( 6\%)\\
\hline
\multirow{4}{*}{\begin{sideways} $\nabla^4$ \end{sideways}}
\T &   200 &      198 ( 0\%)&       57 ( 1\%)&       17 ( 1\%)&       11 ( 3\%)\\
   &   400 &     \NoConverge&      321 ( 0\%)&       86 ( 1\%)&       45 ( 4\%)\\
   &   800 &     \NoConverge&     1565 ( 0\%)&      484 ( 2\%)&      225 ( 5\%)\\
   &  1600 &     \NoConverge&     7900 ( 1\%)&     2171 ( 2\%)&     1012 ( 5\%)\\
\hline
\end{tabular}
\label{table:RuntimeTable}
\end{table}

Table \ref{table:RuntimeTable} shows the wall-clock time spent on setup (excluding partitioning) 
and solution of some problems from table \ref{table:BigResultsTable}.
The 2D problems were solved with a MATLAB implementation that 
solved each subdomain problem with the ``backslash'' operator on each iteration,
but stored a factorization of the coarse grid.
The 3D problems were solved with a parallel C++ implementation
that solved the subdomains with successive over-relaxation,
and solved the coarse grid with conjugate gradient, preconditioned with incomplete Cholesky.
Both implementations ran on a 32-core Intel Xeon E5-2690 with 256GB of RAM.
In all cases, the bulk of the runtime is spent on the subdomain solves,
despite the use of poorly-scaling solvers for the coarse problem.
Furthermore, the runtime scales approximately linearly with problem size, as desired.

\begin{figure}[th]
\centering
\includegraphics[width=0.8\textwidth]{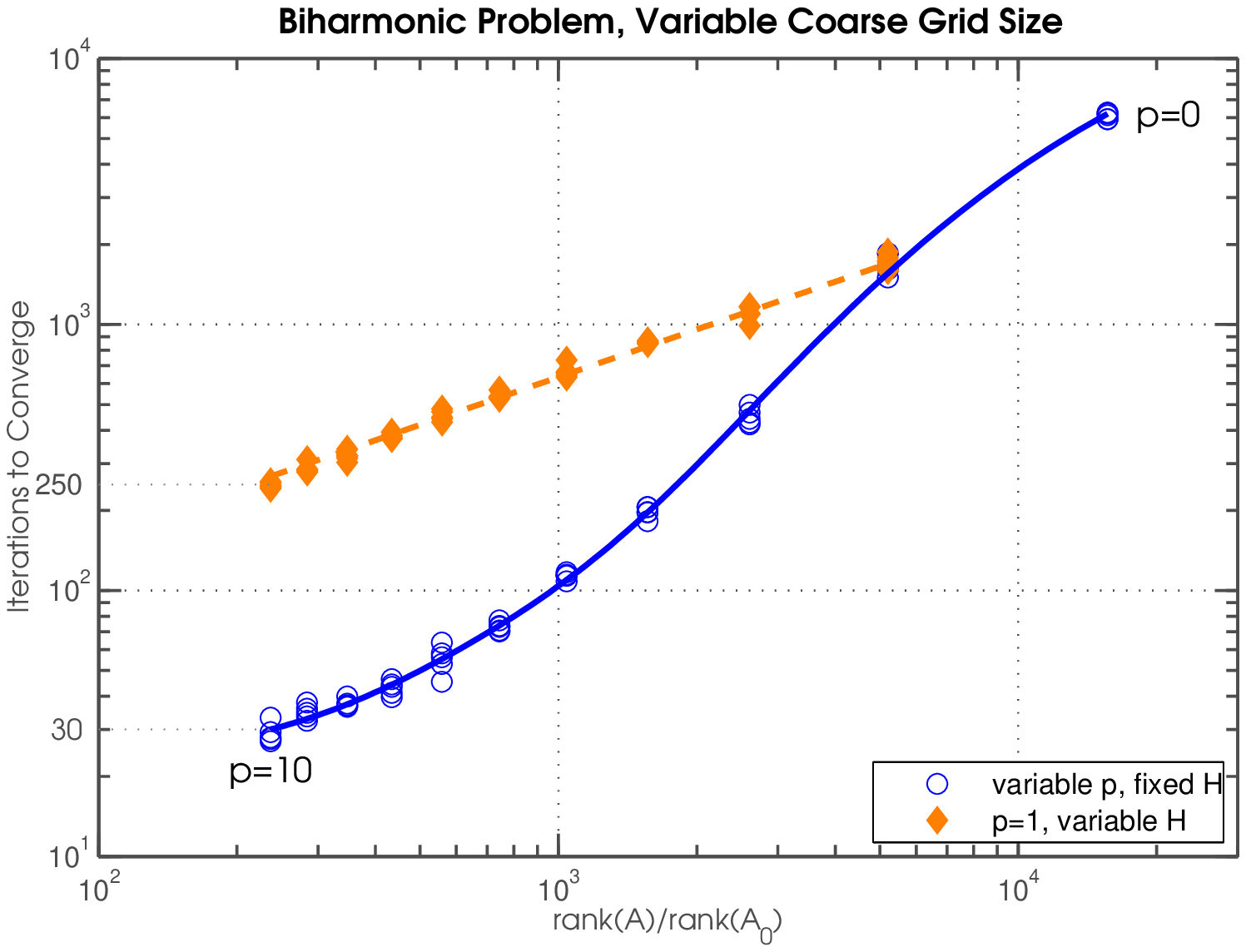}
\caption{
The number of required iterations decreases rapidly with polynomial degree used in the coarse basis 
(parameters: $h^{-1}=1000$, $\Hvh=125$, $\Delta=2$).
For comparison, when using $p=1$ but aggregating using smaller subdomains,
the iteration count is much higher for the same number of coarse variables.
}
\label{fig:p_dependence_biharmonic}
\end{figure}

To directly explore the value in using higher-order coarse bases,
we solve a $1000\times 1000$ biharmonic problem with varying $p$,
a large coarsening factor $\Hvh=125$, and $\Delta=2$.
The large coarsening factor is desirable for efficient parallel implementations,
but significantly reduces the accuracy of the low-order coarse grids.
As shown in figure \ref{fig:p_dependence_biharmonic},
we observe the number of iterations decreases rapidly with increasing $p$.

Note that increasing $p$ increases the size of $\Ao$,
so there are diminishing returns with large $p$.
Nonetheless, significant reductions in problem size are achieved for all $p$.
The least reduction, with $p=10$, is $\rank(\Ao)=(1/284)\cdot\rank(\A)$.
A similar effect occurs in aggregation techniques when the aggregation
subdomains are smaller than the subdomains used in smoothing.
We compare to this approach by using $p=1$ but aggregating on smaller subdomains.
The high-$p$ basis significantly outperforms this approach (figure \ref{fig:p_dependence_biharmonic}).

\begin{figure}[ht]
\centering
\includegraphics[width=0.8\textwidth]{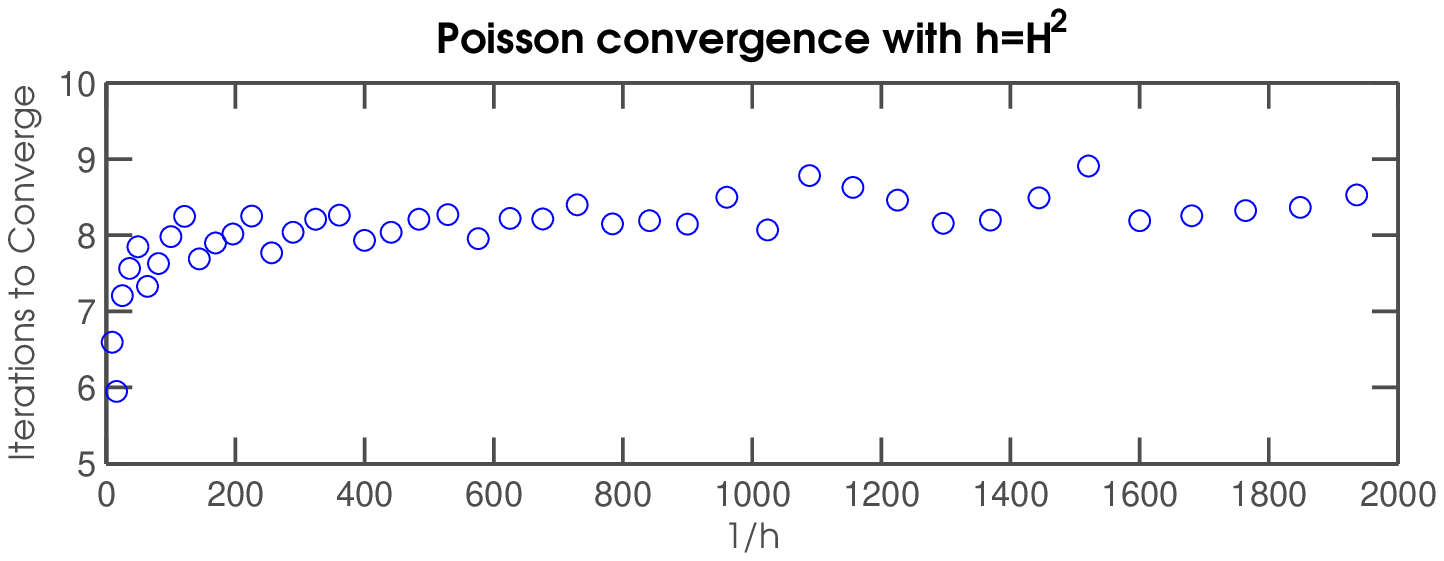}
\caption{
Our grid is convergent and has approximately $h$-independent convergence,
even when $\Hvh$ is not fixed.
Parameters: $h=H^2$, $2\delta\approx H$, $p=1$.
}
\label{fig:h_is_Hsquared}
\end{figure}

Our error bound does not strictly require that $\Hvh=\BigOh(1)$ in order to produce a convergent coarse grid,
and accompanying optimal preconditioner.
For Poisson ($q=1$) using a piecewise linear coarse basis ($p=1$) and substituting the relationship $h=H^2$,
we arrive at the convergent bound
\begin{equation}
\|\Rot\u_0 - \u\|_\A \leq cH^{\sfrac{1}{2}} |u|_{\widetilde{W}_\infty^{2}(\Omega)}.
\end{equation}
Figure \ref{fig:h_is_Hsquared} shows the Poisson problem B with these parameters.
Each subdomain has a number of nodes equal to the number of subdomains,
so the coarse matrix and the subdomain matrices are a similar size,
which is an interesting point in the design space. 
As in all overlapping DD methods, the smoother is very sensitive to $H/\delta$. 
To keep it approximately constant, we set ${\Delta=\lfloor\sfrac{H}{4h}\rfloor}$.
The minor saw-tooth pattern in the graph comes directly from the remaining variation in $H/\delta$.
With higher polynomial degrees, $h=H^s$ for higher powers of $s$ would be possible.

\section{Discussion}

There are a number of outstanding questions raised by this work.
We've shown a bound on the error of the coarse grid,
dependent on $p$ and the smoothness of the solution.
Ideally, we would have a thorough understanding of the 
relationship between all of the parameters ($H$, $h$, $p$, $\delta$, and $q$)
and the condition number or number of iterations to converge.
We leave closing this gap in the analysis for future work.

On the more practical side,
the generating vectors $\F$ are not difficult to supply,
but it would be more convenient to construct similar high-order coarse grids directly from $\A$.
Also, our approach still has the same undesirable dependency on $\Hvh$ that 
is present in many algebraic approach, but is not present in geometric methods.
Following the connection between the coarse basis and DG,
we have done some preliminary work on algebraically constructing a DG-like discretization that is independent of $\Hvh$,
but with mixed success.

\section{Conclusion}
We have presented an algebraic coarse grid construction that
produces a convergent rediscretization of the PDE for a wide variety of PDEs.
It works for both scalar- and vector-valued problems, 
both second- and fourth-order PDEs,
and both smooth and discontinuous coefficients.
The high-order DG-like coarse basis is easy to construct algebraically, 
and Galerkin projection generates a high-order convergent coarse rediscretization of the input problem.
Combined with DD and CG, we observe convergence in a number of iterations nearly independent of problem size.
Furthermore, increasing the polynomial degree rapidly reduces the number of required iterations:
the fastest solves used high-degree polynomials.


\bibliographystyle{siam}
\bibliography{ddg} 

\end{document}